\documentclass[10pt]{article}
\usepackage{oldlfont}
\usepackage{enumerate}
\usepackage{hyperref}
\usepackage{amssymb}
\begin{document}
%%%%%%%%%%%%%%%%%%%%%%%%%%%%%%%%%%%%%%%%%%%%%%%%%%%%%%%%%%%
%  Example 2, Graphs april 8
%%%%%%%%%%%%%%%%%%%%%%%%%%%%%%%%%%%%%%%%%%%%%%%%%%%%%%%%%%%

\pagestyle{headings}                                                    %LMS
\flushbottom                                                            %LMS

\makeatletter                                                           %LMS

% 1. Smaller section titles                                             %LMS
\def\section{\@startsection {section}{1}{\z@}{-1.5ex plus -.5ex         %LMS
minus -.2ex}{1ex plus .2ex}{\large\bf}}                                 %LMS

% 2. Period after a single section number                               %LMS
\renewcommand{\thesection}{\arabic{section}.}                           %LMS
\renewcommand{\thesubsection}{\thesection\arabic{subsection}.}          %LMS
\renewcommand{\thesubsubsection}{\thesubsection\arabic{subsubsection}.} %LMS
\renewcommand{\theparagraph}{\thesubsubsection\arabic{paragraph}.}      %LMS
\renewcommand{\thesubparagraph}{\theparagraph\arabic{subparagraph}.}    %LMS

% -- but avoid a double period                                          %LMS
\def\@thmcountersep{}                                                   %LMS

% 3. Period rather than colon in figure captions                        %LMS
\long\def\@makecaption#1#2{\vskip 10pt \setbox\@tempboxa\hbox{#1. #2}   %LMS
   \ifdim \wd\@tempboxa >\hsize   % IF longer than one line:            %LMS
       #1. #2                  %   THEN set as ordinary paragraph.   %LMS
     \else                        %   ELSE  center.                     %LMS
       \hbox to\hsize{\hfil\box\@tempboxa\hfil}                         %LMS
   \fi}                                                                 %LMS

%4. Put running heads in paper.                                         %LMS
\def\ps@headings{                                                       %LMS
 \def\@oddhead{\footnotesize\rm\hfill\runninghead\hfill}               %LMS
 \def\@evenhead{\@oddhead}                                              %LMS
 \def\@oddfoot{\rm\hfill\thepage\hfill}\def\@evenfoot{\@oddfoot} }      %LMS

%%%%%%%%%%%%%%%%%%%%%%%%%%%%%%%%%%%%%%%
\newcommand{\ZZ}{{\rm Z}\kern-3.8pt {\rm Z} \kern2pt}
\newcommand{\RR}{${\rm I\kern-1.6pt {\rm R}}$}
\newcommand{\NN}{{\rm I\kern-1.6pt {\rm N}}}
\newcommand{\HH}{{\rm I\kern-1.6pt {\rm H}}}
\newcommand{\QQ}{${\rm Q}\kern-3.9pt {\rm \vrule height7.1pt
    width.3pt depth-2pt} \kern7.5pt$}
\newcommand{\CC}{${\rm \kern 3.5pt \vrule height 6.5pt
    width.3pt depth-1.3pt \kern-4pt C \kern.8pt}$}
\newcommand{\FF}{{\rm I\kern-1.6pt {\rm F}}}
\newtheorem{thm}{Theorem}[section]
\newtheorem{lemma}[thm]{Lemma}
\newtheorem{cor}[thm]{Corollary}
\newtheorem{prop}[thm]{Proposition}
\newcommand{\pf}{\noindent {\em Proof.} \ }
\newcommand{\eop}{$_{\Box}$ \vspace{5mm} \relax}
\newtheorem{num}[thm]{}
\newcommand{\split}{\colon\!}
\newcommand{\nsplit}{\cdot}
\newcommand{\Aut}{{\rm Aut}}
\newcommand{\MSC}{{\rm Mathematics Subject Classification: }}
\newcommand{\PSL}{{\rm PSL}}
\newcommand{\PGL}{{\rm PGL}}
\newcommand{\PGaL}{{\rm P\Gamma L}}
\newcommand{\GL}{{\rm GL}}
\newcommand{\SL}{{\rm SL}}
\newcommand{\AGL}{{\rm AGL}}
\newcommand{\ASL}{{\rm ASL}}
\newcommand{\GaL}{{\rm \Gamma L}}
\newcommand{\AGaL}{{\rm A\Gamma L}}
\newcommand{\Sp}{{\rm Sp}}
\newcommand{\ASp}{{\rm ASp}}
\newcommand{\GSp}{{\rm GSp}}
\newcommand{\AGSp}{{\rm AGSp}}
\newcommand{\ZAGSp}{{\rm \hat{A}GSp}}
\newcommand{\nl}{\trianglelefteq}
\newcommand{\mod}{\hbox{mod }}
\newcommand{\Alt}{{\rm Alt}}
\newcommand{\Sym}{{\rm Sym}}
\renewcommand{\o}{\overline}
\newcommand{\sdp}{\rtimes}

\providecommand{\keywords}[1]{\textbf{Keywords:} #1}
\newcommand\blfootnote[1]{%
  \begingroup
  \renewcommand\thefootnote{}\footnotetext{#1}%
  \addtocounter{footnote}{-1}%
  \endgroup}
\newcommand{\Addresses}{{% additional braces for segregating \footnotesize
  \bigskip
  \footnotesize
\noindent
  \textsc{Dipartimento di Matematica e Informatica, Universit\`a della Calabria,
    87036 Arcavacata di Rende, Italy}\par\nopagebreak\noindent
  \textit{E-mail address}: \texttt{jozef.vanbon@unical.it}
}}

\title{A new family of locally $5$-arc transitive graphs of pushing up type with respect to the prime 3}

\author{J. van Bon}

\date{\ }

\maketitle

\Addresses

\def\runninghead{{\sc Locally $5$-arc transitive graphs}}
\pagestyle{headings}

\begin{abstract}
Let $q$ be a power of the prime 3.
A locally 5-arc transitive $G$-graph of pushing up type is constructed for each value of $q$. For $q=3$, the $G$-graph constructed provides an example of a graph with a vertex stabilizer amalgam of shape ${\cal E}_1$ in the sense of \cite{b1}.  Whereas, for the other values of $q$, the vertex stabilizer amalgam of the $G$-graph is of a previously unknown shape. In particular, for $q \neq 3$, these graphs are the first examples of locally 5-arc transitive graphs containing a vertex $z$ for which the group that fixes all 3-arcs originating at $z$ is non-trivial.
\end{abstract}

\blfootnote{
\noindent
Keywords: locally $s$-arc transitive graphs, group amalgams
\par\hskip .35cm
Mathematics Subject Classification (2010): 05C25 \, 05E18 \, 20B25
}

\section{Introduction}

The graphs in this paper will be finite, connected, undirected and without loops or multiple edges.
A $G$-graph is a graph $\Delta$ together 
with a subgroup $G \leq \Aut(\Delta)$. An $s$-arc is a path
$(x_0,x_1,\dots, x_s)$ with $x_{i-1}\neq x_{i+1}$ for $1\leq i \leq s-1$. The stabilizer in $G$ of an $s$-arc $(x_0,x_1,\dots, x_s)$ will be denoted by $G_{x_0,x_1,\dots, x_s}:=\{g \in G \mid x_i^g=x_i, \hbox{ for all } 0 \leq i \leq s \}$.
A $G$-graph $\Delta$ is called 
\begin{itemize}
\item {\em thick} if the valency of each vertex is at least 3;
\item {\em locally finite} if for every vertex $z \in V\Delta$ the stabilizer $G_z$ is a finite group;
\item {\em locally $s$-arc transitive}  if for each vertex  $z\in V\Delta$ the stabilizer $G_z$ acts transitively on the set of arcs of length $s$ originating at $z$.
\end{itemize}

\par\smallskip\noindent
We will assume from now on that $\Delta$ is a locally finite and
locally $s$-arc transitive $G$-graph with $s \geq 1$. Let
$(x_1,x_2)$ be an 1-arc. 
The triple $(G_{x_1}, G_{x_2};G_{x_1,x_2})$ will be called the 
{\em vertex stabilizer amalgam} of $\Delta$, with respect to the 1-arc $(x_1,x_2)$.
Observe that the group $G$ has 
at most 2 orbits on the vertex set $V\Delta$ and is transitive on the edge set $ E\Delta$, since $s \geq 1$.
Hence, the isomorphism type of the vertex stabilizer amalgam does not depend on the choice of the edge. Moreover, we have that
$\Delta$ is either a regular graph, or a bipartite graph with at most 2 valencies.
\par\smallskip 
The vertex stabilizer amalgam 
$(G_{x_1},G_{x_2};G_{x_1,x_2})$ of a $G$-graph $\Delta$ describes the
structure of $G$ and $\Delta$ locally.
Since there is no hope of classifying $G$-graphs from their local properties only, the best one can hope for is to determine the isomorphism type of the vertex stabilizer amalgams possible, and to bound $s$. This is the main problem in the theory of locally $s$-arc transitive graphs. 

This problem finds its origins in the seminal work of  Tutte \cite{T1, T2} on arc-transitive $G$-graphs of valency 3. For arc-transitive graphs of arbitrary valency the problem was extensively studied by Weiss, see \cite{weiss} for an overview,
who also obtained the bound of $s \leq 7$ \cite{weiss8}.
The vertex stabilizer amalgams
for locally finite and locally $s$-arc transitive graphs of valency 3 were determined by Goldschmidt \cite{Gold} (there are 15 isomorphism types, if $s \geq 1$). More
recently it was shown by van Bon and Stellmacher \cite{bonstell} that the vertex stabilizer amalgam for a thick, locally finite and locally $s$-arc transitive $G$-graph is a weak $(B,N)$-pair when $s \geq 6$. Hence, the vertex stabilizer amalgam is known by the classification of weak $(B,N)$-pairs obtained by Delgado and Stellmacher \cite{DS}. These two result combined imply that any thick, locally finite and locally $s$-arc transitive $G$-graph will have $s \leq 9$.

The structure of the vertex stabilizer amalgam of a thick $G$-graph with $s \geq 4$ seems to be restricted too. We need to introduce some more notation and definitions before we can describe what is known in this case.
For a vertex $z \in V\Delta$, the kernel of action of $G_z$ on the set of neighbors $\Delta(z)$ will be denoted by
$G_z^{[1]}$.
A locally finite $G$-graph $\Delta$ is called of
\begin{itemize}
\item {\em local characteristic $p$} with respect to $G$,
if there exists a prime $p$ such that
$$C_{G_{x_i}}(O_p(G_{z}^{[1]})) \leq O_p(G_{z}^{[1]}), \hbox{ for all } z\in V\Delta;$$
\item {\em pushing up type} with respect to the 1-arc $(x_1,x_2)$ and the prime $p$, if $\Delta$ is of local characteristic $p$ with 
respect to $G$, and $O_p(G_{x_1}^{[1]}) \leq O_p(G_{x_2}^{[1]})$.
\end{itemize}

In the first part of the proof of \cite{bonstell} it is established that a thick, locally finite and locally $s$-arc transitive $G$-graph
with $s \geq 6$, has to be of local characteristic $p$, but cannot be of pushing up type.
This is not true any more when $s =5$. However,
the vertex stabilizer amalgams for thick, locally finite and locally $s$-transitive $G$-graphs with $s \geq 4$, that are not of local characteristic $p$ were determined in \cite{TWloc}
Except for 2 vertex stabilizer amalgams, they all belong to one infinite family. In all cases, the corresponding $G$-graphs have $s=5$.
The classification of the vertex stabilizer amalgams of thick, locally finite and locally $s$-arc transitive $G$-graphs of
pushing up type with $s \geq 4$ is in progress, see \cite{bp, b1, b2}. 
We expect that the vertex stabilizer amalgams of thick, locally finite and locally $s$-arc transitive $G$-graphs that are of
local characteristic $p$, but not of pushing up type, 
are weak $(B,N)$-pairs when $s \geq 4$.

There are not many $G$-graphs of pushing up type with $s \geq 4$ known, and those that are known all have $s=5$. By \cite[1.1]{bp}, a $G$-graph of pushing up type will have a vertex $z$ with $|\Delta (z)|=q+1$ and $\PSL_n(q) \nl G_x^{\Delta (z)} \leq \PGaL_2(q)$.
In \cite{b1} it was shown that there are four vertex stabilizer amalgams possible when $q=3$. In \cite{Unitary} two $G$-graphs of pushing up type were constructed which give completions of the amalgams ${\cal D}_2$ and ${\cal E}_2$ of \cite{b1}. In \cite{Symmetric} an infinite family of thick, locally finite and locally $5$-arc transitive $G$-graphs of pushing up type was constructed with $|O_p(G_{x_1,x_2})|=q^3$, with $q \not\equiv 1 (\mod 3)$. The graph for $q=3$ in this family provides a completion of the amalgam of shape ${\cal D}_1$ of \cite{b1}, whereas, when $q > 3$ the corresponding vertex stabilizer amalgams were not known to belong to locally $5$-arc transitive $G$-graphs before.

In this paper we construct for each $q$, where $q$ is a power of the prime 3, a thick, locally finite and locally $5$-arc transitive $G$-graphs of pushing up type with $|O_3(G_{x_1,x_2})|=q^4$.
The construction follows the same line of reasoning as in \cite{Symmetric}. However, this time the vertex stabilizers are constructed as subgroups of $\Sym(q^3)$, and the vertex stabilizers generate a proper subgroup of $\Sym(q^3)$.
The resulting vertex stabilizer amalgams are not isomorphic to those given in \cite{Symmetric} or \cite{Unitary}. For $q=3$, the
vertex stabilizer amalgam provides a completion of an amalgam of shape ${\cal E}_1$ of \cite{b1}.

Let $\FF_q$ be a finite field of order $q$.
The semidirect product of $\SL_2(q)$ and $\GaL_2(q)$ with its natural module are denoted by $\ASL_2(q)$ and $\AGaL_2(q)$, respectively. 
Let $\GaL_2(q,S):=N_{\GaL_2(q)} (S)$, where $S \in Syl_p(\SL_2(q))$.
Furthermore, elementary abelian groups of order $q$ and $q^3$ will be denoted with $E_q$ and $E_{q^3}$, respectively.
Our main result is the following. 

\begin{thm}\label{main1}
Let $q$ be a power of the prime 3. Let $\Delta$ be the graph constructed in section 2. Then $\Delta$ is a locally finite and locally 5-arc transitive $G$-graph of pushing up type, such that, for an edge $\{x_1,x_2\}\in E\Delta$, the following hold:
\begin{enumerate}[(i.)]
\item $|\Delta(x_1)|=q+1$ and $|\Delta(x_2)|= q$;
\item $G_{x_1}^{\Delta (x_1)} \cong \PGaL_2(q)$ and $G_{x_2}^{\Delta (x_2)} \cong \AGaL_1 (q)$;
\item $G_{x_1} \cong E_{q^3} \sdp \GaL_2(q)$;
\item $G_{x_2} \cong (E_q \times U) \sdp \GaL_2(q,S)$, where $U$ is isomorphic to a Sylow $3$-subgroup of $\ASL_2(q)$; 
\item $G_{x_2}$ is transitive on 6-arcs originating at $x_2$.
\end{enumerate}
\end{thm}

\par\noindent
{\bf Remark}. A completion $G$ of the vertex stabilizer amalgam $(G_{x_1},G_{x_2};G_{x_1,x_2})$ is determined in \ref{B4} as a subgroup of $\Sym(q^2) \wr \AGaL_1(q)$. We did not compute the automorphism group of the graph $\Delta$, but it is expected to be $G$. 

\par\smallskip
For certain values of $q$, the vertex stabilizer amalgam of the graph $\Delta$ occurring in \ref{main1} contain a subamalgam which also occur as the vertex stabilizer amalgam of thick, locally finite and locally $5$-arc transitive graphs of pushing up type which are covers of $\Delta$. They are given in the second and third theorem of this paper. Before stating it, we need to introduce some more notation.

Let ${\mathfrak A}: \, H_1\leftarrow H_{1,2} \rightarrow H_2$ be an amalgam.
For $i \in \{1,2\}$, let
$T_i$ be the largest subgroup of $H_{1,2}$ such that $T_i \nl H_i$. 
It was shown in \cite[2.1]{TWloc} that there exists a unique inclusion-minimal subgroup $B$ of $H_{1,2}$
such that 
$$T_1T_2 \leq B \hbox{ and } B=H_{1,2}\cap \langle B^{H_i} \rangle, \hbox{ for } i=1,2.$$
For $i \in \{1,2\}$ we define $H_i^*=\langle B^{H_i} \rangle$. The amalgam ${\mathfrak A}^*: \, H_1^*\leftarrow B \rightarrow H_2^*$  is the called the {\em core} of $\mathfrak A$.

Let $S \in Syl_p(ASL_2(q))$ and ${\cal A}:=\{V \leq S \mid V \hbox{ elementary abelian and }|V|=q\}$.
We define $\ASL_2(q,S):= N_{\ASL_2(q)}(S)$ and 
$\AGL_2(q,S):= N_{\AGL_2(q)}(S)$. Let $K_2(q,S)$ denote the kernel of action of $\AGL_2(q,S)$ on $\cal A$.
Furthermore, let $\GL_2(q,S):=N_{\GL_2(q)}(S)$, where 
$S\in Syl_p(\GL_2(q))$.
For any positive integer $d$ dividing $(q-1)$,
each of the groups $\GL_1(q)$, $\AGL_2(q)$ and $K_2(q,S)$, contain a unique normal subgroup of index $d$.
We will denote these subgroups by $\GL_1(q)^{(d)}$, $\AGL_2(q)^{(d)}$ and $K_2(q,S)^{(d)}$, respectively.

The next two theorems will give the structure of the vertex stabilizer amalgams found in this paper in terms of their core.

\begin{thm}\label{main2}
Let $\Gamma$ be one of the $G$-graphs constructed in Section 4.
Let $(x_1,x_2)$ be an 1-arc of $\Gamma$, with $|\Gamma (x_1)|=q+1$ and $|\Gamma (x_2)|=q$, and let ${\mathfrak A} =(H_{x_1},H_{x_2};H_{x_1,x_2})$ be the vertex stabilizer amalgam of $\Gamma$ with respect to this arc. Moreover, for a vertex $z \in V\Delta$, let $Q_z:=O_3(H_{z}^{[1]})$ and let $L_1:= \langle Q_u \mid u \in \Delta (x_1)  \rangle Q_{x_1}$.
Then there exists a divisor $d$ of $q-1$
such that the following hold for the core ${\mathfrak A}^*$:
\begin{enumerate}[(i.)]
\item $L_1\cong E_q \times ASL_2(q)$, $Z(L_1)\cong E_q$ and $O_3(B)=O_3(L_1 \cap B)=O_3(T_2) \in Syl_3(L_1)$;
\item $B= (L_1 \cap B)T_2$ with $L_1 \cap B \cong E_q \times \ASL_2(q, S)$ and $T_2/Z(L_1) \cong K_2(q,S)^{(d)}$;
\item $H_1^* \cong E_{q^3} \sdp \GL_2(q)^{(d)}$ and $C_{H_1^*}(Z_1) \cong E_q \times \AGL_2(q)^{(\frac{q-1}{2})}$;
\item $H_2^* \cong O_3(T_2) \sdp  (\GL_1(q)^{(d)} \times \AGL_1(q))$ and $C_{H_{x_2}^*}(Z(O_3(T_2))) \cong E_q \times K_2(q,S)^{(\frac{q-1}{2})}$.
Moreover, $H_1/L_1$ is isomorphic to a point stabilizer in a 2-transitive subgroup of $\AGaL_1(q)$.
\end{enumerate}
\end{thm}
\noindent

\par\noindent
{\bf Remark}. The core of the vertex stabilizer amalgam of the graph $\Delta$  of \ref{main1}, appears in \ref{main2} for $d=1$. 

\par\smallskip\noindent
For $q=9$, the vertex stabilizer amalgam of the graph $\Delta$ of \ref{main1}, contains a yet another subamalgam which also gives rise to a locally 5-arc transitive graph. To be more precise, we have the following theorem.

\begin{thm}\label{main3}
There exists a locally 5-arc transitive $H$-graph of pushing up type $\Gamma$, containing an 1-arc $(x_1,x_2)$, with $|\Gamma(x_1)|=10$ and $|\Gamma(x_2)|=9$, such that, for,
the group $L_1:= \langle O_3(H_{u}^{[1]})\mid u \in \Delta (x_1) \rangle O_3(H_{x_1}^{[1]})$,
the vertex stabilizer amalgam ${\mathfrak A} =(H_{x_1},H_{x_2};H_{x_1,x_2})$, and the core ${\mathfrak A}^*$,
the following hold:
\begin{enumerate}[(i.)]
\item $L_1\cong E_3 \times \ASL_2(9)$, $Z(L_1)\cong E_3$ and $O_3(B)=O_3(L_1 \cap B)=O_3(T_2)\in Syl_3(L_1)$;
\item $B =(L_1 \cap B)T_2$ with $L_1 \cap B \cong E_3 \times \ASL_2(9,S)$ and $T_2/Z(L_1) \cong K_2(9,S)^{(2)}$;
\item $H_1^* \cong E_{27} \sdp \GL_2(9)^{(2)}$ and $C_{H_1^*}(Z(L_1)) \cong E_3 \times \AGL_2(9)^{(4)}$; 
\item $H_2^* \cong O_3(T_2) \sdp (\GL_1(9)^{(2)} \times \AGL_1(9))$ and $C_{H_2^*}(Z(O_3(T_2))) \cong E_3 \times K_2(9,S)^{(4)}$;
\item $|H_{x_1}:H_1^*|=|H_{x_2}:H_2^*|=|H_{x_1,x_2}:B|=2$ and $H_{x_1}/L_1 \cong Q_8$.
\end{enumerate}
\end{thm}

\par\smallskip\noindent
{\bf Remark}. Besides the examples given in \cite{Unitary}, \cite{Symmetric}, and the ones constructed here, there is one more $G$-graph of pushing up type with $s \geq 4$ known.
The vertex stabilizer amalgam $(G_{x_1},G_{x_2};G_{x_1,x_2})$ of this graph is realized in $\Aut (He)$, with $G_{x_1}\cong \AGL_2(7)^{(2)}$, $G_{x_2}=N_{\Aut (He)} (S)$, where $S \in Syl_7 (\Aut (He))$, and $G_{x_1,x_2} = N_{G_{x_1}}(S)$, as can be seen from \cite{atlas}.

\par\medskip
The organization of this paper is as follows.
In Section 2 we construct a thick, locally finite and locally 5-arc transitive $G$-graph $\Delta$ as the connected component
of a coset graph on two subgroups of $\Sym(q^3)$.  This will yield Theorem \ref{main1}. In section 3 we determine the structure of the group $G$.
In section 4 we define some covers $\Delta_J$ of $\Delta$, whose connected components are locally 5-arc transitive $H$-graphs, for some 
subgroup $H\leq G$.
Moreover, in this section, Theorem \ref{main2} and Theorem \ref{main3} are proven.
\par\smallskip
The notation used in the paper is standard in the theory of (locally) $s$-arc transitive $G$-graphs.
We refer the reader to \cite{KS} for facts about coprime action, and other results in group theory we use in this paper.

\section{The construction}

Let $\FF_q$ be a field of order $q$, with $q=3^r$ for some integer $r \geq 1$. We first construct two subgroups of $\Sym(q^3)$.
We will then construct the coset graph $\Delta$ on these two subgroups. This graph turns out to be a locally $5$-arc transitive graph of pushing up type. The graph will be bipartite and bivalent, with valencies $q+1$ and $q$.

Let $\Omega =\{(a,b,c) \mid a,b,c \in \FF_q \}$. 
Let $\varepsilon = \frac{q}{3}$. Then $\varepsilon$ is relatively prime with $q-1$, so the map 
$x \mapsto x^{\varepsilon}$ is bijection of $\FF_q$ to itself. 
We define the following permutations on $\Omega$:

\begin{enumerate}[]
\item $\alpha_{(u,v,w)}:\, (a,b,c) \mapsto (a+u,b+v,c+w)$, where $u,v,w \in \FF_q$;
\item $\beta_{\lambda, \mu}:\, (a,b,c) \mapsto (\lambda a, \mu b, (\lambda\mu)^{2 \varepsilon} c)$, where $\lambda, \mu \in \FF_q^*$;
\item $\gamma_e :\, (a,b,c) \mapsto (a+eb,b,c)$, where $e\in \FF_q$;
\item $\delta :\,(a,b,c) \mapsto (b,-a,c)$;
\item $\tau_d :\, (a,b,c) \mapsto (a+db^2+d^{\varepsilon}c,b,c)$, where  $d\in \FF_q$;
\item $\sigma:\, (a,b,c) \mapsto (a^3,b^3,c^3)$.
\end{enumerate}

In the next 3 lemmas we will determine various relations 
between these permutations and use these
to define the subgroups we need to construct the graph.
Permutations will be applied from left to right.

\begin{lemma}\label{A1}
Let $u_1,u_2,v_1,v_2,w_1,w_2,e_1,e_2,d_1,d_2 \in \FF_q$ and 
$\lambda_1, \lambda_2, \mu_1, \mu_2 \in \FF_q^*$.
The following relations hold:
\begin{enumerate}[(i)]
\item $\alpha_{(0,0,0)}=\beta_{1,1}=\gamma_0 =\tau_0 =1$;
\item $\alpha_{(u_1,v_1,w_1)} \alpha_{(u_2,v_2,w_2)} = \alpha_{(u_1+u_2,v_1+v_2,w_1+w_2)}$ and $\alpha_{(u_1,v_1,w_1)}^{-1}= \alpha_{(-u_1,-v_1,-w_1)}$;
\item $\beta_{\lambda_1, \mu_1} \beta_{\lambda_2, \mu_2} = \beta_{\lambda_1\lambda_2, \mu_1\mu_2}$  and $\beta_{\lambda_1, \mu_1}^{-1} = \beta_{\lambda_1^{-1}, \mu_1^{-1}}$;
\item $\gamma_{e_1}\gamma_{e_2}=\gamma_{e_1+e_2}$ and $\gamma_{e_1}^{-1}=\gamma_{-e_1}$;
\item $\delta^{-1}=\delta\beta_{-1,-1}$ and $\delta^2=\beta_{-1,-1}$;
\item $\tau_{d_1}\tau_{d_2} = \tau_{d_1+d_2}$ and  $\tau_{d_1}^{-1}=\tau_{-d_1}$.
\end{enumerate}
\end{lemma}

\pf $(i)$-$(v)$ are straightforward. We will only prove $(vi)$. 
Let $(a,b,c)\in \Omega$. Applying first $\tau_{d_1}$, and then $\tau_{d_2}$, to $(a,b,c)$  we obtain
$$(a,b,c) \mapsto ((a+d_1b^2+d_1^{\varepsilon}c,b,c) \mapsto ((a+d_1b^2+d_1^{\varepsilon}c+d_2b^2+d_2^{\varepsilon}c,b,c) =$$
$$((a+(d_1+d_2)b^2 + (d_1^{\varepsilon}+d_2^{\varepsilon})c,b,c)=
((a+(d_1+d_2)b^2 + (d_1+d_2)^{\varepsilon}c,b,c).$$ 
Hence $\tau_{d_1}\tau_{d_2} = \tau_{d_1+d_2}$, and  $\tau_{d_1}^{-1}=\tau_{-d_1}$ follows.
\eop

We define the following groups:
$$A:= \langle \alpha_{(u,v,w)} \mid u,v,w \in \FF_q \rangle, ~
V:=\langle \alpha_{(u,v,0)} \mid u,v \in \FF_q \rangle, ~
F:=\langle \alpha_{(0,0,w)} \mid w \in \FF_q \rangle, $$

$$Z_0=:\langle \alpha_{(u,0,w)} \mid u,w \in \FF_q \rangle, ~  
Z:=\langle \alpha_{(u,0,0)} \mid u \in \FF_q \rangle, ~ 
C:=\langle \gamma_e \mid e \in\FF_q \rangle$$

$$R:=\langle \beta_{\mu, \mu} \mid \mu \in \FF_q^*\rangle, ~ 
S:=\langle \beta_{\mu^2, \mu} \mid \mu \in \FF_q^*\rangle, ~ 
T:=\langle \beta_{\mu, \mu^{-1}} \mid \mu \in \FF_q^*\rangle, ~
$$

$$D:=\langle \delta \rangle,~ E= \langle \tau_d \mid d \in\FF_q \rangle, ~
\Sigma := \langle \sigma\rangle,$$

$$M:= \langle C,D,T \rangle, ~
Q:=\langle A,C \rangle \hbox{ and } 
P:=\langle Q ,E \rangle.$$

\par\smallskip\noindent
By  \ref{A1} the groups $A$, $V$, $F$ ,$Z_0$, $Z$, $C$ and $E$ are elementary abelian groups of order $q^3$, $q^2$, $q$, $q^2$, $q$, $q$ and $q$, respectively,
$D \cong C_4$ and $R \cong S \cong T \cong C_{q-1}$. Moreover, observe that the group $\Sigma$ is a cyclic group of order $r$.

\begin{lemma}\label{A2}
Let $u,v,e \in \FF_q$ and 
$\lambda, \mu \in \FF_q^*$.
The following relations hold:
\begin{enumerate}[(i)]
\item $\beta_{\lambda, \mu}^{-1}\alpha_{(u,v,w)}\beta_{\lambda , \mu}=\alpha_{(\lambda u,\mu v, (\lambda\mu)^{2 \varepsilon}w)}$;
\item $\beta_{\lambda, \mu}^{-1}\gamma_{e}\beta_{\lambda , \mu} =\gamma_{\lambda\mu^{-1}e}$;
\item $\gamma_e^{-1} \alpha_{(u,v,w)} \gamma_e =  \alpha_{(u+ev,v,w)}$;
\item $\delta^{-1}\alpha_{(u,v,w)} \delta = \alpha_{(v,-u,w)}$;
\item $\delta^{-1} \beta_{\lambda, \mu} \delta = \beta_{\mu, \lambda}$
and $\delta^{-1} \beta_{\mu^2, \mu} \delta = 
\beta_{\mu^{-1},\mu}\beta_{\mu^2, \mu}$.
\end{enumerate}
\end{lemma}

\pf 
The relations follow from an easy calculation using the definitions and is left to the reader \eop

\begin{lemma}\label{A3}
The following hold:
\begin{enumerate}[(i)]
\item $[D,T] \leq T$ and $[C,T]\leq C$;
\item $[V,C]=Z$, $[V,D]=V$ and $[V,T]=V$;
\item $[F,S]\leq F$, $[V,S]\leq V$, $[D,S]=T$, $[T,S]=1$ and $[C,S]=C$;
\item $[V,F]=[D,F]=[T,F]=[C,F]=1$;
\item $[R,C]=[R,D]=[R,S]=[R,T]=1$, $[R,V]=V$ and $[R,A]\leq A$;
\item $[M,F]=1$, $[M,V]=V$ and $[M,S]\leq M$;
\item $MS \cap A=MS\cap V=1$, $C\cap A=1$ and $T \cap S=1$.
\end{enumerate}
\end{lemma}

\pf
$(i)$-$(v)$ follow immediately from  \ref{A2} and coprime action. $(vi)$ follows from $(i)$-$(v)$.
To prove $(vii)$ observe that $A$ acts regularly on $\Omega$, $V \leq A$, and that $C$,$D$,$T$ and $S$ all fix $(0,0,0)$. Hence $MS \cap A=MS\cap V=1$, and so also $C\cap A=1$.
Let $\beta_{\lambda, \mu} \in S \cap T$. Then $\mu^{-1}=\lambda= \mu^2$, hence $\mu^3=1$.
Since $q \not\equiv 1 ~(\mod 3)$, we have that $\mu=1$, and thus 
$S \cap T=id$.
\eop

\begin{lemma}\label{A4}
We have
\begin{enumerate}[(i)]
\item $M \cong SL_2(q)$ and $MS^{(d)} \cong GL_2(q)^{(d)}$;
\item $VM \cong ASL_2(q)$ and $VMS \cong AGL_2(q)$;
\item $\langle A, M, S \rangle \cong (F \times (V \sdp M) )\sdp S \cong (F \times ASL_2(q)) \sdp C_{q-1}$;
\item $\langle A, M, S \rangle \cong (F \times V) \sdp (M \sdp S) \cong (F \times V) \sdp GL_2(q)$;
\item $C_{MS}(F) \cong GL_2(q)^{(\frac{q-1}{2})}$.
\end{enumerate}
$V$ is the natural module for $GL_2(q)$ and $F$ an elementary abelian $p$ group of order $q$.
\end{lemma}

\pf
Observe that, $M$ normalizes $V$ and centralizes $F$, and that $A=VF$.
Note that $A$ is regular on $\Omega$. Thus $C_{Sym(q^3)}(A)$ is regular on $\Omega$ too. 
Hence $|C_{Sym(q^3)}(A)|=|A|$ and, in particular,
$C_{Sym(q^3)}(A)=A$, since $A$ is an abelian group.
From  \ref{A3} it now follows that $C_{M}(A)\leq M \cap A =1$.
Whence $AM \cong A \sdp M$ and $C_M(A)=1$.
Since $[M,F]=1$ we have $C_{M}(V)=C_{M}(VF)=C_{M}(A)=1$.

By  \ref{A1} and  \ref{A2}, 
$R$ centralizes the groups $T$, $C$ and $D$, hence $[R,M]=1$.
Observe also that $M$ normalizes $V$, and that $R$ 
acts fix point freely on $V$.
It follows that $V$ is an $\FF_q M$-module, and since $C_{M}(V)=1$, this action is faithful.

Let 
$e_1=\alpha_{(1,0,0)}$ and $e_2=\alpha_{(0,1,0)}$. Then $(e_1,e_2)$ is an $\FF_q$-basis for $V$. Since $D$, $C$, $T$ induce $\FF_q$-linear transformations on $V$, it follows from  \ref{A2} that all these transformations have determinant 1. Hence $M \cong SL_2(q)$.

Note that $[R,S]=R\cap S=1$. Hence, $V$ is also an $\FF_q MS$-module.
Let $m \in C_{MS}(V)$, then $m=m_0b$, for certain $m_0\in M$ and $b \in S$.
Since both $m$ and $b$ normalize $\langle e_1 \rangle$ and $\langle e_2 \rangle$, we have
that $m_0$ normalizes $\langle e_1 \rangle$ and $\langle e_2 \rangle$ too. Whence $m_0\in T \cap S =1$, by \ref{A3}. 
It follows that $C_{MS}(V) \leq C_S(V)=1$. Note that 
$S$ act as $\FF_q$-linear transformations
on $V$ and $\beta_{\mu^2, \mu}$ has determinant $\mu^3$. Since $q \not\equiv 1 ~(\mod 3)$ it follows that
$M \cap S=1$, $MS\cong GL_2(q)$ and $MS^{(d)}\cong GL_2(q)^{(d)}$.
We conclude that, $\langle A, M, S \rangle \cong (F \times (V \sdp M) )\sdp S \cong (F \times ASL_2(q)) \sdp C_{q-1}$ and
$\langle A, M, S \rangle \cong (F \times V) \sdp MS \cong (F \times V) \sdp GL_2(q)$.

To prove the last statement, observe that $C_{MS}(F)=MC_S(F)$, and $C_S(F)= \langle {\beta_{1,-1}} \rangle$ is a group of order 2.
\eop

We now start constructing the second subgroup we need.
This will involve the groups $P$, $T$ and $S$.
We begin with investigating the groups $Q$ and $E$ and the action of $S$ and $T$ on them.

\begin{lemma}\label{A5}
The following hold:
\begin{enumerate}[(i)]
\item $Q=CA$ and $|Q|=q^4$;
\item $[ \gamma_{e_1}\alpha_{(u_1,v_1,w_1)} , \gamma_{e_2}\alpha_{(u_2,v_2,w_2)}]=\alpha_{(c_2v_1-c_1v_2,0,0)}$, 
where  $u_1,u_2,v_1,v_2,w_1,w_2,e_1,e_2,\in \FF_q$;
\item $[Q,Q]=Z$;
\item $Z(Q)=Z_0$.
\end{enumerate}
\end{lemma}

\pf
$(i)$. By \ref{A3} $A \cap C=1$ and $[A,C]\leq A$, hence $Q=AC$. Since $|A|=q^3$ and $|C|=q$, it follows that  $|Q|=q^4$.
\par\noindent
$(ii)$. This follows from \ref{A2}.
\par\noindent
$(iii)$. This follows from $(ii)$. 
\par\noindent
$(iv)$. By $(i)$ $Q=AC$. 
It follows easily from \ref{A2} that $Z_0 \leq Z(AC)$. On the other hand, if $\alpha\gamma \in Z(AC)$, for some $\gamma \in C$ and $\alpha \in A$, then  $[\langle \gamma \rangle , A]=1$, since $A$ is an abelian group. Hence, $\gamma \in C_C(A)=1$, by \ref{A2} $(iii)$.
It follows that $Z(AC) \leq A$ and from  \ref{A2} $(iii)$ that $Z(AC)=Z_0$.
\eop

\par\noindent
We have shown that $Q$ is a $p$-group of order $q^4$ with $Z(Q)=Z_0$.
Recall that $E$ is an elementary abelian $p$-group of order $q$.
The following lemma lists some relations involving the elements of the group $E$.

\begin{lemma}\label{A6}
Let $u_1,u_2,v_1,v_2,w_1,w_2,e_1,e_2,d \in \FF_q$ and 
$\lambda, \mu \in \FF_q^*$.
The following relations hold:
\begin{enumerate}[(i)]
\item $\beta_{\lambda, \mu}^{-1} \tau_{d} \beta_{\lambda, \mu} = \tau_{\lambda \mu^{-2}d}$;
\item $\tau_d ^{-1} \alpha_{(u_1,v_1,w_1)} \tau_d=  \gamma_{2dv_1}\alpha_{(u_1+dv_1^2+d^{\varepsilon} w_1,v_1,w_1)}$;
\item $\tau_d ^{-1} \gamma_{e_1} \tau_d =\gamma_{e_1}$.
\end{enumerate}
\end{lemma}

\pf
$(ii)$ and $(iii)$ follow from a straightforward calculation using the definitions.
We will only give the proof of $(i)$. Recall that $q-1=3\varepsilon -1$.
Whence, for any $\rho \in \FF_q ^*$, we have
$\rho^{\varepsilon}=\rho^{1-2\varepsilon}$.
%Applying subsequently $\tau_d ^{-1}$, $\alpha_{(u_1,v_1,w_1)}$ and $\tau_d$ to the element $(a,b,c)\in\Omega$ we obtain
Applying subsequently 
$\beta_{\lambda, \mu}^{-1}$, $ \tau_{d}$ and $ \beta_{\lambda, \mu}$ to the element $(a,b,c)\in\Omega$ we obtain

$$
(a,b,c) \mapsto (\lambda^{-1} a, \mu^{-1} b, (\lambda\mu)^{-2\varepsilon} c) \mapsto
(\lambda^{-1} a +d(\mu^{-1}b)^2+ d^{\varepsilon}  (\lambda\mu)^{-2\varepsilon} c, \mu^{-1} b, (\lambda\mu)^{-2\varepsilon} c) \mapsto$$

$$
( a +d\lambda(\mu^{-1}b)^2+ d^{\varepsilon} \lambda (\lambda\mu)^{-2\varepsilon} c, b, c) =
( a +(d\lambda\mu^{-2})b^2+ d^{\varepsilon} \lambda^{1-2\varepsilon} \mu^{-2\varepsilon} c, b, c)=$$

$$
( a +(d\lambda\mu^{-2})b^2+ d^{\varepsilon} \lambda^{\varepsilon}\mu^{-2\varepsilon} c, b, c)=
( a +(d\lambda\mu^{-2})b^2+ (d\lambda\mu^{-2})^{\varepsilon} c, b, c).$$
Whence $\beta_{\lambda, \mu}^{-1} \tau_{d} \beta_{\lambda, \mu} = \tau_{\lambda \mu^{-2}d}$.
\eop

\par\noindent
Next we give the action of $T$ and $S$ on $E$ and $Q$.

\begin{lemma}\label{A7}
Let $u,v,w, e,d \in \FF_q$ and 
$\lambda, \mu \in \FF_q^*$.
The following relations hold:
\begin{enumerate}[(i)]
\item $[\beta_{\mu^2, \mu}, \tau_d] =id$;
\item $\beta_{\mu^2, \mu}^{-1} \gamma_{e}\alpha_{(u,v,w)}\beta_{\mu^2, \mu}= \gamma_{\mu e}\alpha_{(\mu^2 u, \mu v, \mu^2 w)}$;
\item $\beta_{\lambda, \lambda^{-1}}^{-1} \tau_{d} \beta_{\lambda, \lambda^{-1}} = \tau_{\lambda^3 d}$;
\item $\beta_{\lambda, \lambda^{-1}}^{-1} \gamma_{e}\alpha_{(u,v,w)} \beta_{\lambda, \lambda^{-1}}= \gamma_{\lambda^2 e}\alpha_{(\lambda u, \lambda^{-1} v, w)}$.
\end{enumerate}
\end{lemma}

\pf 
$(i)$ and $(iii)$ follow immediately from \ref{A6}. Whereas,
$(ii)$ and $(iv)$ follow from \ref{A2} and the fact that
$6\varepsilon \equiv 2 \, (\mod q-1)$, hence $\mu^{6 \varepsilon} = \mu^2$, for all $\mu\in \FF_q^*$.
\eop

\begin{lemma}\label{A8}
The following hold:
\begin{enumerate}[(i)]
\item $E$ normalizes $Q$ and $E \cap Q=1$;
\item $[T,E]=E$, $T$ is transitive on the non-trivial elements of $E$, and $T$ normalizes $P$;
\item $[F,E]=Z$ and $F$ normalizes $PT$;
\item $[S,E]=1$ and $S$ normalizes $PT$;
\item $SET \cong S \times ET \cong GL_1(q) \times AGL_1(q)$;
\item $PT=(Q\sdp E) \sdp T$;
\item $Z(P)=Z$.
\end{enumerate}
In particular, $|P|=q^5$.
\end{lemma}

\pf
$(i)$.
We deduce from \ref{A6} that $[C,E]=[Z,E]=1$, $[F,E]=Z$ and $[V,E]=ZC$.
Whence $[Q,E] \leq Q$. By  \ref{A2} $A \nl Q$ and
by  \ref{A6} $N_E(A)=1$. Hence $Q\cap E \leq N_E(A)=1$.
\par\noindent
$(ii)$.
By  \ref{A7} $[T,E]=E$ and $T$ is transitive on the non-trivial elements of $E$. Since $T$ normalizes $Q$, $T$ also normalizes $P$.
\par\noindent
$(iii)$. By \ref{A3} $[F,M]=1$, hence
$[F,QT]=1$. From  \ref{A6} we get $[F,E]=Z$.
It follows that $[F,P]=Z$, and so also $[F,PT]=Z$.
\par\noindent
$(iv)$.
It follows from \ref{A7} that $[S,E]=1$. Thus $S$ normalizes $QE$, since $S$ normalizes $Q$ by \ref{A3}.
Since $[S,T]=1$, we conclude that $S$ normalizes $QET=PT$.
\par\noindent
$(v)$. By \ref{A3} and \ref{A8}, $[S,ET]=1$ and $S \cong GL_1(q)$.
By \ref{A8} $ET\cong AGL_1(q)$, so in particular $Z(ET)=1$.
Hence $S \cap ET=1$ and $(v)$ follows.
\par\noindent
$(vi)$. 
Recall that $Q \cap E=1$ by $(i)$.
Since $Q$ is a $p$-group, we have $Q \cap ET \leq Q\cap E=1$, and $(iv)$ follows.
\par\noindent
$(vii)$.
From \ref{A6} we get $[Z,E]=1$, and by \ref{A5}, $Z_0=Z(Q)$.
We conclude that $Z\leq Z(P)$.
Let $\tau \in Z(P)$, then $[\tau, Q]=1$.
Since $Q \nl P$ and $Q/Z$ is an abelian group, we have that $\tau$ induces the identity on $Q/Z$.
We can write  $\tau= \tau_d\gamma_e\alpha_{(u,v,w)}$, 
for some $\tau_d \in E$, $\gamma_e \in C$ and $\alpha_{(u,v,w)}\in A$. 
Since by \ref{A5} $[Q, Q]\leq Z$, $\tau$ and $\tau_d$ induce the same automorphism on $Q/Z$. 
Whence $[\langle \tau_d\rangle , Q]\leq Z$. 
From \ref{A6} we now deduce that $d=0$. 
Whence $\tau \in Q$, and thus $Z(P) \leq Q$.
But then, by \ref{A5} $Z(P) \leq Z(Q)=Z_0=ZF$. 
By $(iii)$ $[E,F]=Z$ and thus $Z(P) \leq Z$. Hence $Z(P)=Z$.
\eop

For any $r\in \FF_q$ we define
$$Q_{r}:=\langle \gamma_{tr}\alpha_{(u,t,w)} \mid t,u,w \in \FF_q\rangle.$$ 
Furthermore, we define
$$Q_*:=CZ_0 \hbox{ and }  \Lambda= \{Q_r \mid r \in \FF_q\}.$$

\par\smallskip\noindent
Note that $Z_0 \leq Q_r$, for all $r\in \FF \cup \{*\}$ and $Q_0=A$.
It follows from 
\ref{A1} and  \ref{A2} that
$Q_*$ and $Q_0$ are elementary abelian $p$-groups of order $q^3$.

\smallskip
The following lemma gives the structure of $Q$ in terms of the groups defined above, and the action of $SET$ on it.

\begin{lemma}\label{A9}
Let $r,s,d\in \FF_q$. The following hold:
\begin{enumerate}[(i)]
\item $T$ is transitive on $\{Q_t \mid t\in \FF_q^*\}$ and normalizes $Q_0$ and $Q_*$;
\item $Q_r$ is elementary abelian of order $q^3$ and $Q_r \nl Q$;
\item $Q_r \cap Q_s= Z_0$ if $r \neq s$ and $Q_r \cap Q_*=Z_0$;
\item $\{Q_t/Z_0 \mid t \in \FF_q \cup \{*\}\}$ is a set of TI-groups of $Q/Z_0$;
\item $[S,Q_*]=Q_*$ and $[S,Q_r]=Q_r$;
\item $\tau_d^{-1} Q_* \tau_d = Q_*$ and $\tau_d^{-1} Q_r \tau_d = Q_{2d+r}$.
\end{enumerate} 
\end{lemma}

\pf
\par\noindent 
\par\noindent
$(i)$. The last two statements follow immediately from  \ref{A7}.
Note that by \ref{A7},
$\beta_{\lambda, \lambda^{-1}}^{-1} \gamma_{tr}\alpha_{(u,t,w)} \beta_{\lambda, \lambda^{-1}}= \gamma_{\lambda^2 tr}\alpha_{(\lambda u, \lambda^{-1} t, w)}
= \gamma_{(\lambda^3r) (\lambda^{-1}t) }\alpha_{(\lambda u, \lambda^{-1}t, w)} \in Q_{\lambda^3r}$. Hence, $\beta_{\lambda, \lambda^{-1}}^{-1} Q_r \beta_{\lambda, \lambda^{-1}} =Q_{\lambda ^3 r}$. 
We conclude that $T$ induces a transitive permutation group on the set $\{Q_r \mid r\in \FF_q^*\}$, since $(3,q-1)=1$.
\par\noindent
$(ii)$. Note that $|Q_2| \geq q^3$, since $Q_2\leq AC$ and $A\cap C=1$.
From \ref{A6} we deduce that $\tau_{-1}^{-1}Q_2\tau_{-1}\leq Q_0$. Whence $Q_2$ is an elementary $p$-group of order $q^3$. 
By $(i)$, $T$ is transitive on the set $\{Q_r \mid r\in \FF_q^*\}$, 
hence $Q_r$ is an elementary $p$ group of order $q^3$, for all $r\in \FF_q$.
Since $Z \leq Q_r$, and by \ref{A5}, $[Q,Q]=Z$, we have $Q_r \nl Q$,
for all $r\in \FF_q$.
\par\noindent 
$(iii)$. By the definition of the group $Q_*$, we have that an element
$\gamma_{c}\alpha_{(u,v,w)} \in Q_*$ if and only if $v=0$.
Whence $Q_* \cap Q_r =Z_0$, for all $r \in \FF_q$.
Since, for all $r\in\FF_q$, $|Q_r|=q^3$, we have that
$Q_r=\{ \gamma_{tr}\alpha_{(u,t,w)} \mid t,u,w \in \FF_q\}$.
Whence $Q_r \cap Q_s=Z_0$ for any $r,s \in \FF_q$ with $r \neq s$.
\par\noindent 
$(iv)$. This follows from $(iii)$.
\par\noindent 
$(v)$. This follows from \ref{A7}.
\par\noindent 
$(vi)$. 
By \ref{A6}, $[E,Q_*]=1$. Let 
$\gamma_{rt}\alpha_{(u,t,w)} \in Q_r \backslash Z_0$.
Then $t\in \FF_q^*$ and by the
aforementioned lemma
$\tau_d^{-1}\gamma_{rt}\alpha_{(u,t,w)}\tau_d=
\gamma_{rt}\gamma_{2dt}\alpha_{(u+dt^2+d^{\varepsilon},t,w)}=
\gamma_{(r+2d)t}\alpha_{(u+dt^2+d^{\varepsilon}w,t,w)} \in Q_{r+2d}$.
Whence $\tau_d^{-1} Q_r\tau_d = Q_{r+2d}$.
\eop

Next we show that the group $\Sigma$ normalizes $F$, $V$, $C$, $D$, $S$, $T$, and $E$. Therefore, it will also normalize
$M$, $MRS$, $Q$ and $PRS$.

\begin{lemma}\label{A10}
Let $u,v,d,e \in \FF_q$ and 
$\lambda, \mu \in \FF_q^*$.
The following relations hold:
\begin{enumerate}[(i)]
\item $\sigma^{-1} \alpha_{(u,v,w)} \sigma = \alpha_{(u^3,v^3,w^3)}$;
\item $\sigma^{-1} \beta_{\lambda, \mu} \sigma = \beta_{\lambda^3, \mu^3}$;
\item $\sigma^{-1} \gamma_{e} \sigma =\gamma_{e^3}$;
\item $\sigma^{-1} \delta \sigma =\delta$;
\item $\sigma^{-1} \tau_{d} \sigma =\tau_{d^3}$.
\end{enumerate}
\end{lemma}

\pf
Straightforward.
\eop

We are now ready to define the graph $\Delta$. Let
$$G_1:=AMS\Sigma, ~ G_2:=QSET\Sigma, G_{1,2}:= G_1 \cap G_2 \hbox{ and } G:=\langle G_1, G_2 \rangle$$
\par\noindent
and let $\Delta$ be the coset graph on $G$-cosets of $G_1$ and $G_2$,
where the two cosets $G_1g_1$ and $G_2g_2$ are called adjacent 
if and only if  $G_1g_1 \cap G_2g_2\neq \emptyset$. 
Let  $x_1=G_1$ and $x_2=G_2$. Observe that $G_{x_1,x_2}=G_{1,2}$.

\begin{lemma}\label{A11} 
The following hold:
\begin{enumerate}[(i)]
\item The action of $G_{x_2}$ on $\Delta (x_2)$ is
equivalent to the action of $ET\Sigma$ on $\Lambda$;
\item$ G_{x_1,x_2}=QST\Sigma$;
\item $G_{x_1}^{[1]}=AR$ and $G_{x_2}^{[1]}=QS$;
\item $G_{x_1}^{\Delta (x_1)}\cong P\Gamma L_2(q)$ and $G_{x_2}^{\Delta (x_2)}=ET\Sigma \cong A\Gamma L_1(q)$;
\item The action of $G_{x_1}^{\Delta (x_1)}$ on $\Delta (x_1)$ is the
natural 2-transitive representation of $P\Gamma L_2(q)$ in characteristic $3$;
\item The action of $G_{x_2}^{\Delta (x_2)}$ on $\Delta (x_2)$ is the
natural 2-transitive representation of $A\Gamma L_1(q)$;
\item $G_{x_1} \cong E_{q^3} \sdp \GaL_2(q)$ and $G_{x_2} \cong (E_q \times U) \sdp \GaL_2(q,S)$, where $U$ is isomorphic to a Sylow $3$-subgroup of $\ASL_2(q)$.
\end{enumerate}
\end{lemma}

\pf
$(i)$ and $(ii)$.
Note that $Q_0=O_3(G_1)$, so $G_{1,2} \leq  N_{G_2}(Q_0)$.
By \ref{A9} and \ref{A10},
we have $N_{G_2}(Q_0)=QST\Sigma \leq G_{1,2}$.
Hence, $G_{1,2} = N_{G_2}(Q_0)=QST\Sigma$ and both $(i)$ and $(ii)$ follow.
\par\noindent
$(iii)$.
By $(ii)$, $QS$ is a subgroup of $G_{1,2}$, and by \ref{A8}, $QS$ is normalized by $ET\Sigma$.
Hence, $QS \leq G_{x_2}^{[1]}$.
By \ref{A8} and \ref{A10}, $G_2/QS \cong ET\Sigma \cong \Gamma L_1(q)$. Since $E$ does
not normalize $Q_0$,  and $G_{x_2}^{[1]}/QS \nl G_2/QS$, we conclude
that $G_{x_2}^{[1]}=QS$.
By \ref{A4} $MS \cong GL_2(q)$. It follows from \ref{A10} that
$MS\Sigma \cong \Gamma L_2(q)$, and so $G_1/A \cong \Gamma L_2(q)$.
By \ref{A3} $M$ normalizes $AR$. Since $AR\leq G_{1,2}$, we have $AR \leq G_{x_1}^{[1]}$.
Since $G_{x_1}^{[1]}/AR \leq G_{1,2}/AR \cong A \Gamma L_1 (q)$, and
$G_{x_1}^{[1]}/AR \nl G_1/AR \cong MS\Sigma/R \cong P\Gamma L_2(q)$ we have
$G_{x_1}^{[1]}=AR$.
\par\noindent
$(iv)$. This follows from $(ii)$, $(iii)$, \ref{A4}, \ref{A8} and \ref{A10}.
\par\noindent
$(v)$ and $(vi)$. These follow from $(i)$ and $(iv)$.
\par\noindent 
$(vii)$ This follows from \ref{A4}, \ref{A8} and \ref{A10}. 
\eop

\begin{lemma}\label{A12} 
$\Delta$ is locally 5-arc transitive $G$-graph of pushing up type.
Furthermore, for any 5-arc $\alpha=(y_0,y_1,y_2,y_3,y_4,y_5)$ we have
\begin{enumerate}[(i.)]
\item $G_{\alpha}/O_3(G_{\alpha})\cong Aut(\FF_q)$;
\item $O_3(G_{\alpha})=Z(O_3(G_{y_2,y_3}))=O_3(G_{y_1}^{[1]})\cap O_3(G_{y_4}^{[1]})$;
\item If $|\Delta(y_0)|=q$, then 
$O_3(G_{\alpha})$ induces a regular group on 
$\Delta (y_5) \setminus \{y_4\}$ with kernel of size $q$. 
In particular, $G_{y_0}$ is transitive on 6-arcs originating at $y_0$. 
\end{enumerate}

\end{lemma}

\pf
Let 
$x_3=G_1 \tau_1$, $x_0=G_2\delta$,
$x_{-1}=G_1\tau_1\delta$ and $x_4=G_2\delta\tau_1$.
Then $x_3 \in\Delta (x_2)$, $x_0\in\Delta(x_1)$, 
$x_{-1}\in \Delta (x_0)$ and $x_{4}\in\Delta (x_3)$,
hence $(x_{-1},x_0,x_1,x_2,x_3,x_4)$ is a 5-arc.
It follows from \ref{A11} that

$$|G_{x_1}: G_{x_1,x_2}|=q+1, \,\, |G_{x_2}: G_{x_1,x_2}|=q \hbox{ and } s \geq 1.$$

\par\noindent
Note that
$G_{x_3}=(G_{x_1})^{\tau_1}$, 
$G_{x_0}=G_{x_2}^\delta$, $G_{{x_{-1}}}=G_{x_3}^\delta$ and 
$G_{x_4}=G_{x_0}^{\tau_1}$.
Moreover, observe that $O_3(G_{x_1})=Q_0$ and $O_3(G_{x_2}^{[1]})=Q$.
Furthermore, we have
$G_{x_3,x_2}=G_{x_1,x_2}^{\tau_1}=\langle F,VC,T^{\tau_1},S\Sigma \rangle$ and
$G_{x_0,x_1}=G_{x_1,x_2}^\delta=\langle A,C^\delta, T, S\Sigma \rangle$.

By \ref{A11} $G_{x_1}^{[1]}=AR$, and 
$G_{x_1}^{\Delta (x_1)}\cong P\Gamma L_2(q)$ is acting in
its natural 2-transitive representation in characteristic $3$ on $\Delta (x_1)$. Therefore, we have that
$G_{x_0,x_1,x_2}=AST\Sigma$. Likewise,
by \ref{A11}, $G_{x_2}^{[1]}=QS$ and $G_{x_2}^{\Delta (x_2)}=ET\Sigma$ is acting in
its natural 2-transitive representation on $\Delta (x_2)$.
Thus, we have that
$G_{x_1,x_2,x_3}=QS\Sigma$.

Hence,

$$|G_{x_1,x_2}:G_{x_1,x_2,x_3}|=q-1 \hbox{ and } 
|G_{x_1,x_2}:G_{x_0,x_1,x_2}|=q.$$

\par\noindent
Since $s \geq 1$, we have that
$G_{x_i}$ is transitive on 2-arcs originating at $x_i$, for $i \in\{1,2\}$.
Whence $s \geq 2$.  

\par\smallskip\noindent
Since $G_{x_0,x_1,x_2,x_3}=G_{x_0,x_1,x_2} \cap G_{x_3}$ normalizes
$Q_2$ and $AS\Sigma \leq G_{x_0,x_1,x_2,x_3}$, it follows that
$G_{x_0,x_1,x_2,x_3}=AS\Sigma$. Hence

$$ |G_{x_0,x_1,x_2}:G_{x_0,x_1,x_2,x_3}|=q-1  \hbox{ and } 
|G_{x_1,x_2,x_3}:G_{x_0,x_1,x_2,x_3}|=q.$$

\par\noindent
Thus, we have that
$G_{x_i}$ is transitive on 3-arcs originating at $x_i$, for $i \in\{0,3\}$.
Whence $s \geq 3$.

\par\smallskip\noindent
We have $G_{x_{-1},x_0,x_1,x_2}=G_{x_0,x_1,x_2,x_3}^\delta =
A(S\Sigma)^\delta$. 
By \ref{A2},
$(S\Sigma)^\delta$ induces $T\Sigma$
on $\Delta (x_2)$. It follows that 
$G_{x_{-1},x_0,x_1,x_2,x_3}=
G_{x_{-1},x_0,x_1,x_2}\cap G_{x_3}= A\Sigma$. 

Since
$G_{x_2,x_3,x_4}=(ATS\Sigma)^{\tau_1}= Q_2ST^{\tau_1}\Sigma$ and
$T^{\tau_1} \leq ET$, we conclude that $T^{\tau_1}$ normalizes $Q_2$ and acts regularly on
$\Delta (x_2) \setminus \{x_3\}$.
Since $Q_2S \leq G_{x_2}^{[1]}$,
it follows that $G_{x_1,x_2,x_3,x_4}=Q_2S\Sigma$. 

Since $Q_2^{\Delta (x_1)} =C^{\Delta (x_1)}$ acts regular on
$\Delta (x_1) \setminus \{x_2\}$, it follows that
$G_{x_0,x_1,x_2,x_3,x_4}= G_{x_0}\cap G_{x_1,x_2,x_3,x_4} =Z_0S\Sigma$.
Hence  

$$|G_{x_0,x_1,x_2,x_3}:G_{x_{-1},x_0,x_1,x_2,x_3}|=q-1 \hbox{ and } 
|G_{x_0,x_1,x_2,x_3}:G_{x_0,x_1,x_2,x_3,x_4}|=q.$$

\par\noindent
Thus, $G_{x_i}$ is transitive on 4-arcs originating at $x_i$, for $i \in\{0,3\}$.
Whence $s \geq 4$.

\par\smallskip\noindent
Since 
$G_{x_{-1},x_0,x_1,x_2,x_3,x_4}=
G_{x_{-1},x_0,x_1,x_2,x_3}\cap G_{x_0,x_1,x_2,x_3,x_4}=
A\Sigma \cap Z_0S\Sigma=Z_0\Sigma$, it follows that

$$|G_{x_0,x_1,x_2,x_3,x_4}:G_{x_{-1},x_0,x_1,x_2,x_3,x_4}|=q-1 
\hbox{ and }
|G_{x_{-1},x_0,x_1,x_2,x_3}:G_{x_{-1},x_0,x_1,x_2,x_3,x_4}|=q.$$

\par\noindent
Thus, $G_{x_i}$ is transitive on 4-arcs originating at $x_i$, for $i \in\{4,-1\}$. Whence $s \geq 5$.

\par\smallskip\noindent
Since $|G_{x_{-1},x_0,x_1,x_2,x_3,x_4}|=q^2|\Sigma|$ and $q_4=q-1>|\Sigma|$, we conclude that $s=5$.
This completes the proof
that $\Delta$ is a locally 5-arc transitive $G$-graph.
By the definitions of the groups $G_1$ and $G_2$ it is clear that it is of pushing up type.
\par\smallskip\noindent
We will now prove $(i) - (iii)$ using the fact that any 5-arc is conjugated to either
$(x_{-1}, x_0,x_1,x_2,x_3,x_4)$ or $(x_4,x_3,x_2,x_1,x_0,x_{-1})$. 
Since $G_{x_{-1},x_0,x_1,x_2,x_3,x_4}=Z_0\Sigma$ and $O_3(G_{x_{-1},x_0,x_1,x_2,x_3,x_4})=Z_0$, $(i)$ follows.
The first equality of $(ii)$ follows from
$O_3(G_{x_{-1},x_0,x_1,x_2,x_3,x_4})= Z_0=Z(O_3(G_{x_1,x_2}))$.
To prove the second equality, observe that
$Z_0=Z(O_3(G_{x_1,x_2}))=Z(O_3(G_{x_2}^{[1]}))$ and $Z_0\leq O_3(G_{x_4}^{[1]})$. Moreover,
$A \cap G_{x_4} = A \cap G_{x_1,x_2,x_3,x_4} = A \cap Q_2S\Sigma =Z_0$. Hence
$A\cap G_{x_4}=Z_0$. It follows that $O_3(G_{x_1}^{[1]})\cap O_3(G_{x_4}^{[1]})=Z(O_3(G_{x_2}))$.
Finally, we prove
$(iii)$.
By \ref{A9} $Q_2 \cap V =Z$, and by \ref{A2} $Z \cap Z^{\delta^{-1}}=1$. 
Since $Z^{\delta^{-1}} \leq V$ it follows that 
$Q_2 \cap Z^{\delta^{-1}}=1$. Hence, $Q_2^{\delta}\cap Z=1$ too,
and thus $Q_0^{\tau_1\delta}\cap Z=1$.
We have shown that $O_3(G_{x_{-1}}) \cap Z=1$. 
We conclude that $O_3(G_{x_0})=O_3(G_{x_{-1}})Z$.
Since $O_3(G_{x_0})$
induces a regular group on
$\Delta (x_{-1}) \setminus \{x_0\}$, we have that $Z$ is regular on
$\Delta (x_{-1}) \setminus \{x_0\}$.
\eop

\par\smallskip\noindent
{\bf Remark}.
Note that, since $F$ fixes the 5-arc $(x_{-1},x_0,x_1,x_2,x_3,x_4)$, and $F \nl G_{x-1}$,
the group $F$ fixes all vertices at distance at most 3 from
$x_1$.

\par\medskip\noindent
{\it Proof of Theorem \ref{main1}}.
The theorem follows from \ref{A11} and \ref{A12}.
\eop

\section{The group $G$}

This section is devoted to determining the structure of the group $G$.

\par\noindent
For $c \in \FF_q$, we define
$$\Omega_c:=\{(a,b,c) \mid a,b \in \FF_q\}.$$
Furthermore, let
$$\widehat \Omega :=\{\Omega_c \mid c \in \FF_q\}, ~ K_1:=VM, ~ K_2:=VCET, ~ 
K_{1,2}:=K_1 \cap K_2 \hbox{ and } 
K:=\langle K_1, K_2 \rangle.$$
\par\noindent

\begin{lemma}\label{B1}
We have
$K_1 \cong ASL_2(q)$, $K_2 \cong Q \sdp AGL_1(q)$ and
$K_{1,2}=VCT$.
\end{lemma}

\pf
Observe that, by \ref{A4} and \ref{A8}, $K_1 \cong ASL_2(q)$ and 
$K_2 \cong Q \sdp AGL_1(q)$.
Clearly $VCT \leq K_{1,2}$. On the other hand,
\ref{A6}  yields $VC \nl K_2$. Hence $VC \nl K_{1,2}$.
Since $CTV/V$ is a maximal subgroup of $K_1/V \cong SL_2(q)$, we conclude that
$VCT=N_{K_1}(VC)\geq K_{1,2}$. Whence $K_{1,2}=VCT$.
\eop

\begin{lemma}\label{B2}
The following hold:
\begin{enumerate}[(i)]
\item $K$ acts trivially on $\widehat \Omega$ and $F$ acts regularly on  $\widehat \Omega$;
\item $FS\Sigma \cong F \sdp (S\Sigma)$, with $S\Sigma \cong \Gamma L_1(q)$ and $C_{S\Sigma}(F)=\langle \beta_{1,-1}\rangle$;
\item the kernel of action of $FS\Sigma$ on $\widehat \Omega$ is $\langle \beta_{1,-1}\rangle$;
\item $FS\Sigma$ normalizes $K_i$ and $K_i \cap FS\Sigma =1$, $i =1,2$.
%$K \cap FS\Sigma \leq  \langle \beta_{1,-1}\rangle$ and
\end{enumerate}
\end{lemma}

\pf 
$(i)$. This follows from
the action of the elements of $V$, $C$, $T$, $D$, $E$ and $F$ on $\widehat \Omega$.
\par\noindent
$(ii)$.
Consider the elements
$\beta_{\mu^2, \mu}\sigma^i \in S\Sigma$ and
$\alpha_{(0,0,f)} \in F$, where $f\in \FF_q$ and $\mu \in \FF_q^*$. 
We have
$(\beta_{\mu^2, \mu}\sigma^i)^{-1}
\alpha_{(0,0,f)}(\beta_{\mu^2, \mu}\sigma^i)=
\alpha_{(0,0,(\mu^6f)^{p^i})}$,
from which we conclude that $C_{S\Sigma}(F)=\langle \beta_{1,-1}\rangle$.
Since $F$ is an abelian group, we have $F\cap S\Sigma \leq C_{S\Sigma}(F) \cap F=1$. From \ref{A10} we deduce $S\Sigma \cong \PGaL_1(q)$.
\par\noindent
$(iii)$.
For $f\in \FF_q$, $\mu \in \FF_q^*$ and $1\leq i \leq r$ we have that
$\alpha_{(0,0,f)}\beta_{\mu^2, \mu}\sigma^i$ maps
$\Omega_c$ to $\Omega_{(f+c\mu^2)^{3^i}}$.
Therefore,
$\alpha_{(0,0,f)}\beta_{\mu^2, \mu}\sigma^i$
acts trivially on $\widehat \Omega$ if and only if 
$f=0$, $i=0$ and $\mu \in \{1,-1\}$.
\par\noindent
$(iv)$.
By  \ref{A10} the group $\Sigma$ normalizes
$F$, $V$, $C$, $D$, $T$, $S$, and $E$.
By \ref{A3} $F$ and $S$ normalize
$V$, $C$, $D$, $T$, and $E$.
Hence, $FS\Sigma$ normalizes both $K_1$ and $K_2$.
The group $K$ acts trivially on $\widehat \Omega$, by $(i)$ .
Hence, by $(iv)$ $K \cap FS\Sigma \leq \langle \beta_{1,-1} \rangle$.
By the action of $\beta_{1,-1}$  on $V$, $\beta_{1,-1} \not\in K_1$.
Since $\beta_{1,-1}$ normalizes $Q_0$ and, by \ref{A7},
$N_{K_2} (Q_0)=QN_{ET}(Q_0)=K_{1,2}$,
it follows that $\beta_{1,-1} \not\in K_2$. 
Hence, $K_i \cap FS\Sigma=1$ for $i\in\{1,2\}$.
\eop

The structure of $K_{1,2}$ and $K_i$, $i=1,2$, was given in \ref{B1}.
The structure of $G_{1,2}$ and $G_i$, $i=1,2$, was given in \ref{A11}.
We determine the structure of $K$ and $G$. 

\begin{lemma}\label{B3} 
The following hold:
\begin{enumerate}[(i)]
\item $G_1 \cong (F \times K_1) \sdp S\Sigma$;
\item $G_2 \cong (K_2 \sdp F) \sdp S\Sigma$ with $[F,K_2]=Z$;
\item $G_{1,2} \cong (F \times K_{1,2})\sdp S\Sigma$;
\item $K \nl G$ and $G=KFS\Sigma$.
\end{enumerate}
\end{lemma}

\pf
$(i)$.
By \ref{A3} $[F, K_1]=1$, and by \ref{B2} $F\Cap K_1 \le F \cap K=1$.
Note that $S\Sigma$ and $K_1$ both stabilize $\Omega_0$, and 
$F$ acts regularly on $\widehat \Omega$. Hence
$(F \times K_1) \cap S\Sigma = K_1 \cap S\Sigma =1$,
by \ref{B2}.
\par\noindent
$(ii)$.
Recall that $Z_0=FZ$.
By \ref{A3} and \ref{A8}, $[F,K_2]=Z \leq K_2$. Hence, $F$ normalizes $K_2$.
Moreover, by \ref{B2} $F \cap K_2 \leq F \cap K=1$.
By \ref{B2},
$K_2$ and $F$ act on $\widehat \Omega$  trivially and regularly, respectively.
Since
$S\Sigma$ stabilizes $\Omega_0$ and $K_2 \cap S\Sigma  =1$, by \ref{B2}.
We obtain $K_2F \cap S\Sigma \leq K_2 \cap S\Sigma  =1$ too.
\par\noindent
$(iii)$. By \ref{A11}  $G_{1,2}=QST\Sigma$, and by \ref{B1} $K_{1,2}=VCT$. Hence,
$G_{1,2}=FK_{1,2}S\Sigma$. The claim now follows from $(i)$ and $(ii)$. 
\par\noindent
$(iv)$.
By \ref{B2} $FS\Sigma$ normalizes $K_1$ and $K_2$, 
hence $K$ is normalized by $FS\Sigma$ and $G=KFS\Sigma$.
\eop

\begin{lemma}\label{B4}
The following hold:
\begin{enumerate}[(i)]
\item $K\langle \beta_{1,-1} \rangle$ is the kernel of action of $G$ on $\widehat \Omega$;
\item $K\cong \Alt(q^2) \times \cdots \times \Alt (q^2)$ ($q$ terms) and $C_G(K)=1$;
\item $G \cong K \sdp (F \sdp C_{\frac{1}{2}(q-1)} \sdp C_r)$ if $q \equiv 1 ~ (\mod ~4)$;
\item $G \cong (K \sdp C_2) \sdp (F \sdp C_{\frac{1}{2}(q-1)} \sdp C_r)$ if $q \equiv 3 ~(\mod ~4)$.
\end{enumerate}
\end{lemma}

\pf
Let $X$ be the kernel of action of $G$ on $\widehat \Omega$.
$(i)$. 
By $(i)$ $G=KFS\Sigma$ and 
by \ref{B2} $K\leq X$ and $F$ act regularly on $\widehat \Omega$.
Hence, $G \cap X =K(FS\Sigma \cap X) =K\langle \beta_{1,-1} \rangle$ by \ref{B2}.
\par\noindent
$(ii)$.
For a subgroup $Y \leq G$ stabilizing $\Omega_0$, we denote by $\overline Y$ 
the group induced by $Y$ on $\Omega_0$. By \ref{B2} $K$ stabilizes $\Omega_0$.
Since $V$ acts regularly on $\Omega_0$, the group $VM$ induces 
a group isomorphic to $ASL_2(q)$ on $\Omega_0$,
which acts 2-transitively on the elements of $\Omega_0$.
Moreover, it is easy to see that $\overline E$ is of order $q$.
Since $E$ does not normalize $V$, we have that 
$\overline {\langle VM , E \rangle}$ is a 
2-transitive permutation group
on $q^2$ elements not normalizing $V$. 
By the classification of finite 2-transitive groups, see \cite{cameron} for a list, $\overline K$ is either
isomorphic to $Alt(q^2)$ or $Sym(q^2)$.
Since $V$, $C$ and $E$ have odd order, and $\overline T$ and $\overline D$ 
are generated by an even permutation on $\Omega _0$,  it follows that
$\overline K \cong \Alt(q^2)$.

Note that $S\Sigma$ fixes $\Omega_0$ and has 2 orbits of length $\frac{1}{2}(q-1)$ on $\widehat \Omega$. Since $F$ acts regularly on $\widehat \omega$, we conclude that $FS\Sigma$ acts primitively on $\widehat \Omega$.
Let $Y$ be a minimal normal subgroup of $G$.
Suppose  $[Y,K]=1$, then $Y \leq C_G(K) \leq C_G(A)= A \leq K$,
for $A$ in an abelian group which acts regularly on $\Omega$ ,
which is a contradiction.
Hence $[Y,K]\neq 1$, and $Y \leq K$. Since $F$ normalizes $Y$ and is regular
on $|\widehat \Omega|$, we must have $\overline Y \neq 1$.
In particular $1 \neq \overline Y \nl \overline K$, and so
$\overline Y \cong \Alt(q^2)$. Consequently
 $\overline {C_G(Y)} =1$. Since $C_G(Y) \nl G$ and $F$ act transitive on $\widehat \Omega$, we must have $C_G(Y)=1$.
Since $Y$ is isomorphic to the direct product of isomorphic simple groups, it follows that
$Y$ is isomorphic to the direct product of groups isomorphic to $\Alt(q^2)$ and
$G$ 

Let $Y_0$ be a component of $Y$ with $Y_0^{\Omega_0} \neq 1$.
Let $Y_i \in E(Y)$ with $Y_i \neq Y_0$. Since $[Y_i,Y_0]=1$ we have
$Y_i^{\Omega_0}=1$. Since $C_G(Y)=1$ we must have that $G$ acts transitive on the components of $E(Y)$. In particular,
$\{\Omega_c \mid Y_0^{\Omega_c} \neq 1\}$ is a block of imprimitivity
for the action of $FS\Sigma$. Hence, the set $\{\Omega_c \mid Y_0^{\Omega_c} \neq 1\}$
is either $\widehat\Omega$ or $\{\Omega_0\}$.
In the first case $Y=E(Y)=Y_0\cong \Alt (q^2)$. Since $C_G(Y)=1$, we must have
that $F$ is isomorphic to a subgroup of $\Aut(Y)$. This is impossible, since $F \cap X=1$, $Y \leq X$ and $q \geq 3$.
Therefore, we are in the latter case.
So $Y_0^{\Omega_c}=1$, for all $c \in \FF_q^*$. It follows that,
$Y \cong Y_0 \times \cdots \times Y_q$, with $Y_i^{\Omega_c}=1$, if and only if $i \neq c$ and $Y_i \cong \Alt (q^2)$.
Since $\overline K= \overline Y_0$ and $G$ is faithful on $\Omega$, we have $K=Y$. 

\par\noindent
$(iii)-(iv)$.
From $(i)$  and \ref{B2} we have that
$X= K \langle\beta_{1,-1} \rangle$
and $G/X \cong F \sdp C_{\frac{1}{2}(q-1)}\sdp C_r$.
Now consider the permutation $\beta_{1,-1}$. 
On each of the sets $\Omega_c$ it has $q$ fixed points. 
Hence, $\beta_{1,-1}$ acts on $\Omega_c$ as the product of $\frac{1}{2}q(q-1)$ 2-cycles.
Whence, induces an even permutation on $\Omega_c$ if $q \equiv \,1\, (\mod ~4)$,
and an odd permutation on $\Omega_c$ if $q \equiv \,3\, (\mod ~4)$.
Hence, $ \beta_{1,-1} \in K$ if $q \equiv \,1\, (\mod ~4)$, and 
$ \beta_{1,-1} \not\in K$ if $q \equiv \,3\, (\mod ~4)$.
Therefore, $X=K$ if $q \equiv \,1\, (\mod ~4)$, and 
$X \cong K \sdp C_2$ if $q \equiv \,3\, (\mod ~4)$.

Assume first that $q \equiv \,1\, (\mod ~4)$.
Let $\zeta$ be a generator of $\FF_q^*$, and write $q-1=2^lm$ with $m$ odd.
Consider the permutation 
$$\theta: (a,b,c) \mapsto \left\{
\begin{array}{rc}
(\zeta^2a, \zeta b, c)& {\rm if}~b \neq 0\\
(\zeta^ma, 0, c) & \rm{otherwise}\\
\end{array} 
\right.$$
Then $\theta$ is of order $q-1$ and $\theta^{\frac{q-1}{2}}=\beta_{1,-1}$. 
Moreover, on each of the sets $\Omega_c$ it induces
$q$ cycles of length $q-1$ and $m$ cycles of length $2^l$, hence it is an even permutation, and thus $\theta \in K$.
Observe that $\theta$ commutes with $\beta_{\zeta^2, \zeta}$ and normalizes $F$.
Now $\theta\beta_{\zeta^2, \zeta} \in G$ has order $\frac{q-1}{2}$,
$\langle F, \theta\beta_{\zeta^2, \zeta}  \rangle \cap X=1$ and 
$(iii)$ follows.
Now assume that $q \equiv \,3\, (\mod ~4)$. Then 
$S= \langle \beta_{1,-1} \rangle \times \langle \beta_{\zeta^4,\zeta^2} \rangle$.
Hence $F\langle \beta_{\zeta^4,\zeta^2} \rangle \cap X=1$  and $(iv)$ follows. 
\eop

\section{Some locally 5-arc transitive covers}

We retain the notation introduced in the previous section.
We define some subgroups $K \leq G(\widehat J) \leq G$ which
still act locally 5-arc transitive on $\Delta$. These subgroups will give rise to the various amalgams mentioned in the second main theorem.

We have $CS\Sigma \cong A \Gamma L_1(q)$. 
Consider $A \Gamma L_1$ in its 2-transitive permutation representation on $q$ points.
Let $C \leq \widehat J \leq CS\Sigma$  such that
$\widehat J$ is a 2-transitive group subgroup of $\AGaL_1(q)$. 
Let $J$ be the stabilizer in $\widehat J$ of a $\gamma_0$, so $J \leq S\Sigma$.
Let $d=|S:S\cap J|$. Since $S$ is cyclic of order $q-1$, we have that $|S \cap J|$ has order $\frac{q-1}{d}$.
Let $\zeta$ be a generator of $S$, then 
$ \langle \zeta^d \rangle = S \cap J$.
We denote is group by $S^{(d)}$.

Observe that, by \ref{B2}, $F$ and $J$ both normalize $K_1$ and $K_2$
and, by \ref{A10}, that $J$ normalizes $F$.

We define

$$G_1(\widehat J):=K_1FJ, ~ G_2(\widehat J):=K_2FJ, ~  G_{1,2}(\widehat J):=G_1(\widehat J) \cap G_2(\widehat J) \hbox{ and } G(\widehat J):=\langle G_1(\widehat J), G_2(\widehat J) \rangle.$$

\par\smallskip\noindent
Note that
$G(\widehat J)=KFJ$. Moreover, we have
$G_1=G_1(S\Sigma)$, $G_2=G_2(S\Sigma)$ and $G=G(S\Sigma)$.

\begin{lemma}\label{D1}
 $G_{1,2}(\widehat J)=K_{1,2}FJ$.
\end{lemma}

\pf
We have $K_1FJ \cap K_2FJ=(K_1\cap K_2FJ)FJ$.
Clearly $K_{1,2}FJ \leq G_1(\widehat J) \cap G_2(\widehat J)$.
By \ref{A6} and \ref{A10} $Q_0 \nl  G_2(\widehat J)$. In particular,
$G_1(\widehat J) \cap G_2(\widehat J)$ normalizes $Q_0$.
By \ref{A9}, $E$ is regular on $\{Q_r \mid r\in\FF_q\}$, hence 
$G_1(\widehat J) \cap G_2(\widehat J) \leq Q_0TS\Sigma$. Therefore,
$G_1(\widehat J) \cap G_2(\widehat J)\leq Q_0TS\Sigma \cap K_2FJ=
K_{1,2}FJ$.
\eop

\begin{lemma} \label{D2}
The following hold:
\begin{enumerate}[(i)]
\item $\beta_{1,-1} \in S^{(d)}$; 
\item $K \cap S\Sigma = K\cap J$;
\item $G(\widehat J)/K \cong FJ/(K\cap J)$.
\end{enumerate}
\end{lemma}

\pf
$(i)$.
Let $q=3^r$ and  write $r=2^ab$, with $b$ odd.
Consider the ring $\ZZ / 2^{a+1} \ZZ$.
Since $3^b$ is an odd number, $[3^b]_{2^{a+1}} \in (\ZZ / 2^{a+1} \ZZ)^*$. Since the 
latter is a group of order $2^a$, we have $(3^b)^{2^a} \equiv 1 \hbox{ mod } 2^{a+1}$.
Hence, $2^{a+1}$ divides $q-1$.
In particular, $2^{a+1}$ divides $|J|$. Since
$|J|/|S^{(d)}|$ divides $f$, it follows that $|S^{(d)}|$ is
of even order, hence $\beta_{1,-1} \in S^{(d)}$, 
for $S^{(d)}$ is cyclic and $\beta_{1,-1} \in S$ has order 2.
\par\noindent
$(ii)$.
By \ref{B2} and $(i)$ we have $S\Sigma \cap K \leq \langle \beta_{1,-1} \rangle\leq J$. Hence, $S\Sigma \cap K= J \cap K$.
\par\noindent
$(iii)$.
Since $K_1$ and $K_2$ are normalized by $F$, $S$ and $\Sigma$, we have
$K \nl G(\widehat J)$ and $G(\widehat J)/K \cong FJ/K$. Since, by \ref{B2}, $K \cap FJ= K \cap J$, the claim follows.
\eop

\begin{lemma}\label{D3}
Let $i \in \{1,2\}$. Then
\begin{enumerate}[(i)]
\item $G(\widehat J) \cap G_i= G_i(\widehat J)$ and $|G_i:G_i(\widehat J)|=|S\Sigma: J|$;
\item $|G:G(\widehat J)|=|S\Sigma: J|$.
\end{enumerate}
\end{lemma}

\pf
$(i)$.
Clearly  $G_i(\widehat J) \leq G(\widehat J) \cap G_i$, for $i \in \{1,2\} $. On the other hand, let $g \in G(\widehat J) \cap G_i$.
By \ref{D2}, $G(\widehat J)/K\cong KFJ/KF$. Thus,
$gKF=jKF$ for some $j \in J$.
Since $J \leq G_i$ and, by \ref{B2}, $KF \cap FS\Sigma =F(K \cap FS\Sigma)=F(K\cap S\Sigma)=F(K \cap J)$, we have that 
$$j^{-1}g \in KF \cap G_i=KF\cap (K_iFS\Sigma) =
K_iF(KF \cap S\Sigma)=K_iF(K \cap J)\leq K_iFJ=G_i(\widehat J).$$
Therefore, $g \in G_i(\widehat J)$, for $J\leq G_i(\widehat J)$. It follows that $G(\widehat J) \cap G_i= G_i(\widehat J)$.
Since $G_i=K_iFS\Sigma$ and $G_i(\widehat J)=K_iFJ$, we obtain $|G_i:G_i(\widehat J)|=|S\Sigma: J|$.
\par\noindent
$(ii)$.
By \ref{D2} and \ref{B2}, $G(\widehat J)/K \cong FJ/(K\cap J)$
and
$G/K \cong KFS\Sigma /K \cong  FS\Sigma/(FS\Sigma \cap K) \cong
FS\Sigma/(K \cap S\Sigma) \cong FS\Sigma/(K \cap J) $. Therefore, we have that
$|G:G(\widehat J)|=|S\Sigma: J|$.
\eop

\begin{lemma}\label{D4}
Let $x=G_1$, $y=G_2$ and $u=G_2\delta$.
The following statements are equivalent:
\begin{enumerate}[(i)]
\item $J$ is transitive on $\Delta (u) \setminus \{x\}$;
\item $CJ$ is a 2-transitive group on $\Delta (x) \setminus \{y\}$;
\item $|J:J \cap \Sigma|=q-1$.
\end{enumerate}
\end{lemma}

\pf

For any $c,d \in \FF_q^*$ let $y_c =G_2\delta\gamma_c$ and $u_d =G_1\tau_d\delta$.
Since  $C$ acts transitively on $\Delta (x) \backslash \{y\}$, $E$ acts transitively on $\Delta (y)$, and
$\delta$ maps $y$ to $u$, it follows that
$$\Delta (x)=\{u,y\} \cup \{ y_c \mid c \in \FF_q^*\} \hbox{ and }\Delta (u)=\{x\} \cup \{ u_d \mid d \in \FF_q^*\}.$$
The group $J$ stabilizes the arc $(u,x,y)$ and acts on both of the above sets.
Let $j \in J$ and $\mu \in \FF_q^*$ and $0 \leq i \leq r$ such that
$j=\beta_{\mu^2,\mu} \sigma^i$. Then
$$y_cj = G_2\delta\gamma_c \beta_{\mu^2,\mu} \sigma^i=G_2\delta\gamma_{(c\mu)^{3^i}}= y_{(c\mu)^{3^i}}$$
and
$$u_dj= G_1\tau_d\delta\beta_{\mu^2,\mu} \sigma^i=G_1\tau_{(d\mu^3)^{3^i}}\delta = u_{(d\mu^3)^{3^i}}.$$

Since $q \equiv 0 \,(\mod 3)$,
the map $\phi: \Delta (x) \setminus \{u,y\} \rightarrow  \Delta (u) \setminus \{x\}$ 
given by $\phi: \,y_l \mapsto u_{l^3}$ is an isomorphism of sets which commutes with the action of $J$.
Hence, $J$ acts transitively on the set $\Delta (u) \setminus \{x\}$ if and only if
$J$ acts transitively on the set $\Delta (x) \setminus \{u,y\}$.
That is, if and only if $CJ$ is a 2-transitive group on $\Delta (x) \setminus \{y\}$.
This proves the equivalence of $(i)$ and $(ii)$.

Moreover, $j$ centralizes $\gamma_1$ if and only if $\mu=1$. Hence $C_J(\gamma_1) = J \cap \Sigma$.
Since $CJ$ is a 2-transitive group if and only if the $J$-orbit of $\gamma_1$ has length $q-1$,
the equivalence of $(ii)$ and $(iii)$ follows.
\eop

Let $\Delta_J$  be the coset graph on $G(\widehat J)$-cosets of $G_1(\widehat J)$ and $G_2(\widehat J)$, where the two cosets $G_1(\widehat J)g_1$ and $G_2(\widehat J)g_2$ are called adjacent if and only if  $G_1(\widehat J)g_1 \cap G_2(\widehat J)g_2\neq \emptyset$.

\begin{lemma} \label{D5}
$\Delta_J$ is a connected locally 5-arc transitive $G(\widehat J)$-graph.
Moreover, if $|\Delta_J(x_0)|=q$ then $G(\widehat J)_{x_0}$ is transitive on 6-arcs originating at $x_0$. 
\end{lemma}

\pf
The graph $\Delta_J$ is connected, since $G(\widehat J)=KFJ=\langle K_1, K_2\rangle FJ=\langle K_1FJ, K_2FJ\rangle =\langle G_1(\widehat J),G_2(\widehat J) \rangle$ .
By \ref{D1} and \ref{D3}, $(G_1(\widehat J),G_2(\widehat J); G_{1,2}(\widehat J))$ is a subamalgam of index
$|S\Sigma :J|$ in $(G_1,G_2; G_{1,2})$. 
By \cite[2.1]{Symmetric} $\Delta_J$ is a $|S\Sigma:J|$-cover of $\Delta$ and
the projection map is a local isomorphism of $G$-graphs.
Let $\beta =(y_0,y_1,\dots, y_5)$ be the
5-arc of $\Delta$ given by
$$(G_1\tau_1\delta, G_2\delta,G_1, G_2, G_1 \tau_1,G_2\delta\tau_1).$$
By \ref{A12},
$G_{\beta}/O_3(G_{\beta})$ is isomorphic to a subgroup of $\Sigma$ and
$O_3(G_{\beta})=Z(O_3(G_{1,2}))$.
Since $\Sigma \leq G_{\beta}$ we have $G_{\beta}=ZF\Sigma$.
Let $\alpha = (x_0,x_1,\dots, x_5)$ be a 5-arc of $\Delta_J$
with $x_2=K_1FJ$ and $x_3=K_2FJ$ such that
$\tau(\alpha)=\beta$.
Then $G_{\alpha}=K_{1,2}FJ\cap G_{\beta}=VCTFJ \cap ZF\Sigma=ZF(\Sigma \cap J)$.
By \ref{D4} $|\Sigma : \Sigma \cap J|=|S\Sigma:J|$.
Hence,
$|G_{\beta}:G(J)_{\beta}|= |\Sigma : \Sigma \cap J|=|S\Sigma : J|=|G:G(J)|$.
The lemma now follows from  \cite[2.2]{Symmetric} and \ref{A12}.
\eop

We determine the shape of the 
vertex stabilizer amalgams of the $G(\widehat J)$-graphs $\Delta_J$.
Let $H_1=G_1(\widehat J)$, $H_2=G_2(\widehat J)$ and $H_{1,2}=G_{1,2}(\widehat J)$. 

\begin{lemma}\label{D6}
$H_1^*=FK_1S^{(d)}$, $H_2^*=QS^{(d)}ET$ and $B=QTS^{(d)}$.
\end{lemma}

\pf
By \ref{A11}, $G_x^{[1]} = AR$ and $G_y^{[1]} = QS$.
Note that $Q=O_3(H_2)=O_3(G_y^{[1]})$.
Since $T \cap S\Sigma=1$, we have $TJ \cap S=J \cap S= S^{(d)}$.
Hence, $H_2^{[1]}=H_{1,2} \cap G_y^{[1]}=Q(TJ \cap S)=QS^{(d)}$. 
Therefore $QS^{(d)} \leq B$.
Since $K_1^{\Delta (x)} \cong PSL_2(q)$ and $QS^{(d)} \leq B$
we have that $K_1FS^{(d)} \leq H_1^*$. Thus $QTS^{(d)} \leq B$.
On the other hand, $ \langle (QTS^{(d)})^{H_1} \rangle=FK_1S^{(d)}$ and
$H_{1,2} \cap FK_1S^{(d)}=QTS^{(d)}$. Hence $B=QTS^{(d)}$, and since
$S \cap T=1$, it follows that $H_1^*=FK_1S^{(d)}$.
It now easily follows that $H_2^*=QS^{(d)}ET$.
\eop

\begin{lemma}\label{D7}
The following hold
\begin{enumerate}[(i.)]
\item $L_1\cong E_q \times ASL_2(q)$, $Z_1\cong E_q$ and $O_3(B)=O_3(L_1 \cap B)=O_3(T_2) \in Syl_3(L_1)$;
\item $B= (L_1 \cap B)T_2$ with $L_1 \cap B \cong E_q \times \ASL_2(q, S)$ and $T_2/Z_1 \cong K_2(q,S)^{(d)}$;
\item $H_1^* \cong E_{q^3} \sdp \GL_2(q)^{(d)}$ and $C_{H_1^*}(Z_1) \cong E_q \times \AGL_2(q)^{(\frac{q-1}{2})}$;
\item $H_2^* \cong O_3(T_2) \sdp  (\GL_1(q)^{(d)} \times \AGL_1(q))$ and $C_{H_2^*}(Z(O_3(T_2))) \cong E_q \times K_2(q,S)^{(\frac{q-1}{2})}$.
\end{enumerate}
\end{lemma}

\pf
$(i)$.  Since  $Q \leq L_1$  we have $FAM\leq L_1$. On the other hand,
$Q\in Syl_3(FAM)$ and $FAM \nl G_1$. Hence $L_1 =FAM$.
The rest follows from \ref{A4} and \ref{D6}.
\par\noindent
$(ii)$. By \ref{D6}, $B=QTS^{(d)}$. It follows from \ref{A11} that
$T_2=QS\cap B=QS^{(d)}$. Since $L_1=FAM$ and $FVM\cap B=QT$ we have
$(L_1\cap B)T_2=QTS^{(d)}=B$. Since $O^p(L_1)=VM$ it
follows that $O^3(L_1) \cap B= VCT \cong ASL_2(q,U)$
\par\noindent
$(iii)$. The first statement follows from \ref{D6} and \ref{B2}, the last
from \ref{A4} and \ref{A10}.
\par\noindent
$(iv)$. This follows from \ref{D6} and \ref{B2}.
\eop

It follows from \ref{B4} that the above amalgam $(H_1,H_2;H_{1,2})$ has a completion in a proper subgroup of $\Sym(q^2) \wr \AGL_1(q)$.

\par\smallskip
For $q=9$ the vertex stabilizer amalgam for the $G(J)$-graph with $J < \Gamma L_1 (9)$  and $J \not\leq GL_1(9)$ contains a subamalgam of index
3 which will also give rise to a locally 5-arc transitive graph. 
Observe first that
$\Gamma L_2(9) / SL_2(9) \cong C_{8} \sdp C_2$ .
Let $SL_2(9) \leq N \leq \Gamma L_2(9)$ be a subgroup of index 2 in $\Gamma L_2(9)$.
Then $N/SL_2(9) \in \{C_8, D_8, Q_8\}$.
Let $\Gamma \hat L_2(9)$ denote the group with $\Gamma \hat L_2(9) / SL_2(9) \cong Q_8$.
The group $\Gamma \hat L_2(9)$ has three subgroups of index 2, one of which is $GL_2(9)^{(2)}$. The two 
other subgroups are conjugated in $\Gamma L_2(9)$ and their isomorphism type we denote by  $\Gamma \hat L_2(9)^{(2)}$.
Let $L_1\cong C_3 \times ASL_2(9)$ and let ${\mathfrak F}=(H_1,H_2;H_{1,2})$ denote the amalgam with
\begin{enumerate}[]
\item $H_1\cong C_3 \sdp  A\Gamma \hat L_2(9)$ and $C_{H_1}(Z(L_1)) \cong C_3 \times A\Gamma \hat L_2(9)^{(2)}$; 
\item $H_2\cong U \sdp (C_4 \times AGL_1(9))\cdot C_2$, with $U \in Syl_3(L_1)$;
\item $H_{1,2} \cong U\sdp (C_4 \times GL_1(9))\cdot C_2$.
\end{enumerate}

\begin{lemma}\label{D8}
The amalgam ${\mathfrak F}$ is the vertex stabilizer amalgam of a locally 5-arc transitive $G$-graph with valencies $10$ and $9$.
\end{lemma}

\pf
Assume $q=9$. Let $\zeta$ be a generator of $\FF_9$.
Let $\beta:=\beta_{\zeta^2, \zeta}$. Then $\beta \in S\Sigma$ and 
$\beta^4=\beta_{1,-1}=(\beta\sigma)^2$. Let
$J:=\langle \beta^2,\beta\sigma \rangle$, then $J \cong Q_8$. It follows from \ref{A2} and \ref{D1} that $J$ acts regular on $C$.
By \ref{D5}, $(K_1FJ, K_2FJ;K_{1,2}FJ)$ is the vertex stabilizer amalgam of the locally 5-arc transitive $G(J)$-graph $\Delta_J$ for $q=9$.

Let $F_0=\langle \alpha_{(0,0, \zeta)} \rangle$, which is a group of order 3. By \ref{A2} and \ref{D1} $\beta^2$ inverts $F_0$, and $\beta\sigma$ centralizes $F_0$, hence $J$ normalizes $F_0$.
By \ref{A3} and \ref{A8}, $[K_1, F_0]=1$ and $[K_2,F_0]\leq Z$. It follows that
$|K_1FJ:K_1F_0J|=3=|K_2FJ:K_2F_0J|$. We have $K_1FJ\cap K_2FJ=(K_1\cap K_2FJ)FJ=(K_1\cap K_2)FJ=K_{1,2}F_0J$
and $|K_{1,2}FJ:K_{1,2}F_0J|=3$.
Since $K_1F_0J \cap K_{1,2}FJ=K_{1,2}F_0J=K_2F_0J \cap K_{1,2}FJ$ the amalgam
$(K_1F_0J, K_2F_0J;K_{1,2}F_0J)$ is a subamalgam of index 3 in 
$(K_1FJ, K_2FJ;K_{1,2}FJ)$.
Since the stabilizer of a 5-arc with middle edge $\{K_1(J), K_2(J)\}$ in $G(J)$ is equal to $FZ$, we deduce that
that $(K_1F_0J, K_2F_0J;K_{1,2}F_0J)$ is the vertex stabilizer amalgam of a locally 5-arc transitive graph too.

Let $H_1=K_1F_0J$. Then,
$H_1\cong C_3 \sdp \Gamma \hat L_2(9)$ and $C_{H_1}(Z(L_1)) \cong C_3 \times \Gamma \hat L_2(9)^{(2)}$,
$H_2\cong U \sdp (C_4 \times AGL_1(9))\cdot C_2$ and
$H_{1,2} \cong U \sdp (C_4 \times GL_1(9))\cdot C_2$.
\eop

\noindent
{\it Proof of Theorem \ref{main2}}.
This is \ref{D7}.

\medskip\noindent
{\it Proof of Theorem \ref{main3}}.
This follows from \ref{D8}.

\vspace{1cm}

\end{document}